\documentclass[a4paper]{article}
\usepackage{amsmath,amssymb,euscript}
\usepackage[latin1]{inputenc}
\usepackage[dvips]{epsfig}

\newcommand{\RR}{\mathbb{R}}

\newcommand{\QQ}{\mathbb{Q}}
\newcommand{\Aff}{\mathsf{A}}
\newcommand{\PP}{\mathcal{P}}
\newcommand{\GG}{\mathsf{G}}
\newcommand{\FF}{\mathsf{F}}
\newcommand{\sym}{\mathfrak{S}}
\newcommand{\cur}{\leftarrow}
\newcommand{\prelie}{\operatorname{PreLie}}
\newcommand{\assoc}{\operatorname{As}}
\newcommand{\ndom}{\operatorname{End}}
\newcommand{\quot}{\operatorname{Mu}}
\newcommand{\vect}{\operatorname{Vect}_\QQ}

\newtheorem{prop}{Proposition}
\newtheorem{lemma}{Lemma}
\newtheorem{defi}{Definition}
\newenvironment{proof}{\begin{trivlist}\item{\bf{Proof.}}}
  {\hfill\rule{2mm}{2mm}\end{trivlist}}
                                                                            
\title{Rooted trees and an exponential-like series} \date{\today}
\author{Frédéric Chapoton}

\begin{document}
\maketitle
\begin{abstract}
  This paper deals with a group of generalized power series associated
  to any augmented operad, focusing on the case of the $\prelie$
  operad. The solution of flow equations using the pre-Lie structure
  on vector fields on an affine space gives rise to an interesting
  element of this group.
\end{abstract}



\setcounter{section}{-1}

\section{Introduction}

This paper is composed of three main ingredients: the first is the
definition of a group of power series $\GG_\PP$ for any augmented
operad $\PP$. This group acts on all $\PP$-algebras having a suitable
completeness property. The Lie algebra of this group, or rather a
graded version of it, has been introduced before, from another point
of view, by Kapranov and Manin \cite{kapr-manin}.

The rest of the paper focuses on the case of the $\prelie$ operad
\cite{rooted}. The group $\GG_{\prelie}$ is a group of series indexed
by rooted trees. The second component of the paper is given by two
quotient operads of $\prelie$ corresponding respectively to linear
trees and corollas. The operad on linear trees is no other than the
associative operad $\assoc$, whereas the operad $\quot$ on corollas
seems not to have been considered before. By functoriality of the
group construction, one gets group maps $\GG_{\prelie}\to
\GG_{\assoc}$ and $\GG_{\prelie}\to \GG_{\quot}$. Whereas
$\GG_{\assoc}$ is isomorphic to the group of formal power series for
the composition product, $\GG_{\quot}$ is related to the group of
formal power series for the pointwise multiplication product.

The final ingredient comes from flow equations for vector fields on
the affine space $\Aff^n$. We recall the formula which gives the
Taylor expansion of the solution of a flow equation using the pre-Lie
structure on the vector space of vector fields on $\Aff^n$. We have
not traced the precise origin of this formula, although it is
certainly not new, see \cite{brouder} and references therein. This
formula can be interpreted as a distinguished element of the group
$\GG_{\prelie}$, which we choose to denote by $\exp^*$, for its image
in $\GG_{\assoc}$ is $\exp-1$, where $\exp$ is the usual exponential
function.

The plan of the paper is as follows: in the first section, we define
the functor $\GG$ from augmented operads to groups. The second section
is about vector fields on affine space and the pre-Lie formula for
flows. The third and fourth sections are devoted to the two quotient
operads of $\prelie$ and the images of $\exp^*$ and its inverse in the
associated groups. Then we give the first terms of the expansion of
$\exp^*$ and its inverse.

\section{A group associated to an augmented operad}

Let $\PP$ be an operad in the category $\vect$ of vector spaces over
$\QQ$ and assume that $\PP(0)=\{0\}$ and that $\PP(1)=\QQ e$
where $e$ is the unit of $\PP$. Such an operad is called augmented.

Let $\FF\PP=\oplus_n \PP(n)_{\sym_n}$ be the direct sum of the
coinvariant spaces, which can be identified with the underlying vector
space of the free $\PP$-algebra on a generator $v$, and
$\widehat{\PP}={\prod}_n \PP(n)_{\sym_n}$ be its completion.

We now define a product on $\widehat{\PP}$. Let $x=\sum_m x_m$,
$y=\sum_n y_n$ be two elements of $\widehat{\PP}$. Choose any
representatives $\bar{x}_m$ of $x_m$ (resp.  $\bar{y}_n$ of $y_n$) in
the operad $\PP$. Then one defines
\begin{equation}
  \label{product}
  x \times y = \sum_{m\geq 1} \, \sum_{n_1,\dots,n_m \geq 1} 
  \langle \gamma(\bar{x}_m,\bar{y}_{n_1},\dots,\bar{y}_{n_m})\rangle,
\end{equation}
where $\langle\,\rangle$ is the quotient map to the coinvariants and
$\gamma$ is the composition map of the operad $\PP$. This product is
indeed well defined, as can be directly checked by using the
equivariance axiom of operads to compare two different choices of
representatives.

\begin{prop}
  The product $\times$ defines the structure of an associative monoid
  on the vector space $\widehat{\PP}$. Furthermore, this product is
  $\QQ$-linear on its left argument.
\end{prop}
\begin{proof}
  Let us first prove the associativity. On the one hand, one has
  \begin{multline}
    (x\times y)\times z=\sum_m \sum_{p_1,\dots,p_m}
    \langle\gamma(\overline{(x\times y)}_m,\bar{z}_{p_1}, \dots,
    \bar{z}_{p_m})\rangle\\=
    \sum_m\sum_{n_1,\dots,n_m}\,\sum_{p_1,\dots,p_{n_1+\dots+n_m}}
    \langle\gamma(\gamma(\bar{x}_m,\bar{y}_{n_1},\dots,
    \bar{y}_{n_m}),\bar{z}_{p_1},
    \dots, \bar{z}_{p_{n_1+\dots+n_m}})\rangle.
  \end{multline}

  On the other hand, one has 
  \begin{multline}
    x\times (y \times
    z)=\sum_m\sum_{n_1,\dots,n_m}\langle\gamma(\bar{x}_m,\overline{(y\times
      z)}_{n_1},\dots,\overline{(y\times z)}_{n_m})\rangle
    \\=\sum_m\sum_{n_1,\dots,n_m}\sum_{(q_{i,j})}\langle\gamma(\bar{x}_m,\gamma(\bar{y}_{n_1},\bar{z}_{q_{1,1}},\dots,\bar{z}_{q_{1,n_1}}),\dots,\gamma(\bar{y}_{n_m},\bar{z}_{q_{m,1}},\dots,\bar{z}_{q_{m,n_m}}))\rangle.
  \end{multline}
  Using then the ``associativity'' of the operad, one gets the
  associativity of $\times$.
  
  It is easy to check that the image of the unit $e$ of the operad
  $\PP$ is a two-sided unit for the $\times$ product.
  
  The left $\QQ$-linearity is clear from the formula (\ref{product}).
\end{proof}

On can characterize invertible elements in this monoid. 
\begin{prop}
  An element $y$ of $\widehat{\PP}$ is invertible for $\times$ if and
  only if the first component $y_1$ of $y$ is non-zero.
\end{prop}
\begin{proof}
  The direct implication is trivial. The reverse one is proved by a
  very standard recursive argument.
\end{proof}

Let us call $\GG_{\PP}$ the set of invertible elements of
$\widehat{\PP}$ for the $\times$ product.

\begin{prop}
  $\GG$ is a functor from the category of augmented operads to the
  category of groups.
\end{prop}
\begin{proof}
  The functoriality follows from inspection of the definitions of
  $\widehat{\PP}$ and $\times$.
\end{proof}

In fact, one can see $\GG_\PP$ as the group of $\QQ$-points of a
pro-algebraic group. The Lie algebra of this pro-algebraic group is
given by the usual linearization process on the tangent space
$\widehat{\PP}$, resulting in the formula
\begin{equation}
    \label{liebracket}
  [x, y] = \sum_{m\geq 1} \sum_{n \geq 1} 
  \langle \bar{x}_{m} \circ \bar{y}_{n}
    -\bar{y}_{n} \circ \bar{x}_{m}\rangle,
\end{equation}
where
\begin{equation}
  \bar{x}_m \circ \bar{y}_n=\sum_{i=1}^{m} \gamma(\bar{x}_m,
\underbrace{e,\dots,e}_{i-1\text{ units}},\bar{y}_n,e,\dots,e).
\end{equation}

The graded Lie algebra structure on $\FF\PP$ defined by the same
formulas has already appeared in the work of Kapranov and Manin on the
category of right modules over an operad \cite[Th.
1.7.3]{kapr-manin}. They explained that it acts by polynomial vector
fields on the underlying vector space of any $\PP$-algebra. For the
endomorphism operad $\ndom(L)$ of a vector space $L$, $\GG_{\ndom(L)}$
is the group of formal diffeomorphisms of $L$ preserving the origin.
So $\GG_\PP$ acts by formal diffeomorphisms on any $\PP$-algebra $L$.

In order to get an (non-formal) action of $\GG_\PP$, one can assume
that the $\PP$-algebra $L$ is complete filtered, \textit{i.e.} endowed
with a decreasing filtration $L=F^1 L\supset F^2 L \supset \dots$
which is complete and compatible with the $\PP$-action.

\section{The flow of a vector field}

Let $\Aff^n$ be the affine space of dimension $n$ over $\RR$. It is
well-know \cite{matsu,note-cras} that there is a structure of pre-Lie
algebra on the vector space $\mathsf{V}_n$ of smooth vector fields on
$\Aff^n$. More precisely, let $x_1,\dots,x_n$ be coordinates on
$\Aff^n$. Given two vector fields $F=\sum F_i \partial_i$ and $F'=\sum
F'_j \partial_j$, their pre-Lie product $F \cur F'$ is given by
\begin{equation}
  \label{def-coord}
  F \cur F' := \sum\sum F'_j ( \partial_j F_i ) \partial_i.
\end{equation}
This does not depend on the choice of affine coordinates.

Let $F\in\mathsf{V}_n$ be a vector field. The flow equation of $F$ is
the following equation for a smooth function $g : \Aff^1 \to \Aff^n$ of the
variable $t$ :
\begin{equation}
  \label{flow}
  \begin{cases}
    \frac{dg}{dt}=F(g),\\
    g(0)=g_0,
  \end{cases}
\end{equation}
where $g_0$ is any chosen point in $\Aff^n$. The smoothness of $F$
ensures unicity of the solution.

One can give a formal Taylor development at $t=0$ of the solution of
(\ref{flow}) using the following construction in the pre-Lie algebra
of vector fields (see \cite[Prop. 4]{note-cras}). Let $F^{\cur k}$
denote the $k^{th}$ right iterate of $F$ for the pre-Lie product,
\textit{i.e.}
\begin{equation}
  F^{\cur k}=\underbrace{((\cdots((F \cur F)\cur F)\cdots \cur F) \cur F)}_{k \text{ factors}}.
\end{equation}

\begin{prop}
  The solution of the flow equation (\ref{flow}) has the following
  formal Taylor expansion :
  \begin{equation}
    \label{taylor}
    g(t)=g_0+\sum_{k\geq 1} F^{\cur k}(g_0) \frac{t^k}{k!}.
  \end{equation}
\end{prop}

The proof is by iterated use of the following lemma, which in turn is
easily proved with the definition (\ref{def-coord}) of the pre-Lie
product.
\begin{lemma}
  Let $g$ be the solution of (\ref{flow}) and $G$ be a vector field on
  $\Aff^n$. Then
  \begin{equation}
    \frac{d G(g)}{dt}=(G\cur F)(g).
  \end{equation}
\end{lemma}

In the special case where the vector field is linear, \textit{i.e.}
$F=\sum a_{i,j} x_i\partial_j$ with $a_{i,j}\in\RR$, the formula
(\ref{taylor}) gives back the classical exponential formula :
\begin{equation}
  g(t)=g_0+(\exp(t A)-I) g_0=\exp(t A) g_0.
\end{equation}
where $A$ is the matrix $(a_{i,j})$.

One can define an element $\exp^* (v)$ of the completed free pre-Lie
algebra on a generator $v$ over $\QQ$ by
\begin{equation}
  \label{recursive}
  \exp^* (v):=\sum_{k\geq 1} \frac{v^{\cur k}}{k!}.
\end{equation}

\begin{defi}
  The image of $\exp^*(v)$ by the usual identification between the
  completed free pre-Lie algebra on $\{v\}$ and $\widehat{\prelie}$ is
  an element of the group $\GG_{\prelie}$, which is denoted by
  $\exp^*$. Its inverse is denoted by $\log^*$.
\end{defi}

The group $\GG_{\prelie}$ acts on the completed free pre-Lie algebra
on $v$, and the action of $\exp^*$ on $v$ is $\exp^*(v)$.

If $F$ is a vector field on $\Aff^n$, then the difference between the
flow at time $0$ and the flow at time $1$ can be considered as another
vector field $G$ on $\Aff^n$. In fact, formula (\ref{taylor}) says
that it is formally given by $\exp^* F$. So the meaning of $\log^*$ is
the reverse operation : knowing the displacement $G$ between time $0$
and time $1$, $\log^* G$ formally recovers the vector field $F$. For
these statements to make sense, one must work within a complete
pre-Lie algebra of vector fields, for example the pre-Lie algebra
$\prod_{k \geq 2} \QQ x^k \partial_x$.

\section{Linear trees and composition}

As shown in \cite{rooted}, the $\prelie$ operad can be described in
terms of labeled rooted trees. By convention, edges are oriented
towards the root. Let us call linear the trees that do no branch, that
is to say all vertices have at most one incoming edge.

We recall here briefly (see \cite{rooted} for more details) the
definition of the composition of two labeled rooted trees $T$ and $S$
on the vertex sets $I$ and $J$ respectively. Let $i\in I$; the
composition of $S$ at the vertex $i$ of $T$ is given by
\begin{equation}
  T\circ_i S=\sum_f T\circ_i ^f  S,
\end{equation}
where the sum runs over all maps $f$ from the set of incoming edges of
the vertex $i$ of $T$ to the set of vertices of $S$, and $T\circ_i ^f
S$ can be described as follows: replace the vertex $i$ by the tree
$S$, grafting back the subtrees of $T$ previously attached to $i$,
according to the map $f$.

\begin{prop}
  The subspace of $\prelie$ spanned by non-linear labeled trees is an ideal.
  The quotient map $\phi$ coincides with the usual map from $\prelie$
  to the associative operad $\assoc$.
\end{prop}
\begin{proof}
  Using the description above of the composition map of the operad
  $\prelie$, it is clear that the composition of two labeled trees, at
  least one of which is non-linear, is again non-linear. The quotient
  operad, spanned by labeled linear trees, has dimension $n!$ in rank
  $n$. Its composition is given by insertion, and can be easily
  identified with the associative operad $\assoc$. The quotient map is
  then checked on generators of $\prelie$ to be the same as the usual
  map.
\end{proof}

By functoriality of the group construction, there is a map, still
denoted by $\phi$, from $\GG_{\prelie}$ to $\GG_{\assoc}$.

\begin{prop}
  The group $\GG_{\assoc}$ is isomorphic to the group of invertible
  formal power series in $x\QQ[[x]]$ for the composition product.
\end{prop}
\begin{proof}
  It is more convenient here to work at the monoid level with
  $\widehat{\assoc}$ and $x\QQ[[x]]$. The vector space
  $\assoc(n)_{\sym_n}$ is one dimensional for all $n$, with a basis
  given by the image by $\phi$ of the linear tree with $n$ nodes. Let
  us denote this basis element by $\theta_n$. By left linearity of
  both monoids, it is sufficient to check the product rule for
  $\theta_m$ and $f=\sum_{n\geq 1} f_n \theta_n$. One finds that
  \begin{equation}
    \theta_m \times f= \sum_{n_1,\dots,n_m \geq 1} f_{n_1}\dots f_{n_m} 
    \theta_{n_1+\dots+n_m},
  \end{equation}
  which proves that the linear map defined by $x^n \mapsto \theta_n$
  is an isomorphism between the monoids $\widehat{\assoc}$ and
  $x\QQ[[x]]$. The proposition follows by taking invertible elements.
\end{proof}

\begin{prop}
  The image of $\exp^*$ by $\phi$ is $\exp x-1$ and the image of
  $\log^*$ is $\log(1+x)$.
\end{prop}
\begin{proof}
  One has to prove that the coefficient of the linear tree with $n$
  nodes in $\exp^*$ is $1/n!$ for all $n$. This can be done by recursion,
  using (\ref{recursive}).
\end{proof}

This proposition explains what happens for linear vector fields. The
linearity of a vector field $F$ on $\Aff^n$ implies that the pre-Lie
algebra generated by $F$ inside $\mathsf{V}_n$ is associative, so that
$\exp^*$ reduces to the standard exponential minus one.

\section{Corollas and pointwise multiplication}

There is another interesting kind of trees, opposite to linear trees.
Let us call corollas the trees of depth no greater than two, where the
depth is the maximum number of vertices in a chain of adjacent
vertices starting from the root.

\begin{prop}
  The subspace of $\prelie$ spanned by labeled non-corollas is an
  ideal.
\end{prop}
\begin{proof}
  Using the description above of the composition of $\prelie$, one
  shows that the depth of the composition of two labeled trees is
  greater or equal than the maximum of the depths of these labeled
  trees. Therefore, if one of the labeled trees has depth greater or
  equal to three, so does the composition.
\end{proof}

One can give a simple description of the quotient operad $\quot$. It
has dimension $n$ in rank $n$ with basis given by the image of the
labeled corollas with $n$ nodes. Let us call $\mu^{n}_i$ the image of
the corolla with $n$ nodes and with root labeled by $i$ for
$i=1,\dots,n$.

Then $\mu_{1}^{1}$ is the unit of $\quot$ and the composition is given
by
\begin{equation}
  \begin{cases}
    \mu^{n}_i \circ_i \mu^{\ell}_j= \mu^{n+\ell-1}_{i+j-1},\\
    \mu^{n}_i \circ_h \mu^{\ell}_j= 0 \text{ for } h\not=i\text{ and
    }\ell\geq 2.
  \end{cases}
\end{equation}

Let $G_1$ be the group of formal power series of the form
$1+x\QQ[[x]]$ for the pointwise multiplication product and $G_2$ be
the multiplicative group $\QQ^*$.  There is an action of $G_2$ on
$G_1$ by substitution: $\lambda \cdot f(x)=f(\lambda x)$. A group
similar to the semi-direct product group $G_2 \ltimes G_1$ has been
considered in \cite[\S 2]{frab-brou}.

From the description of $\quot$ above, one deduces that
\begin{prop}
  The group $\GG_{\quot}$ is isomorphic to $G_2 \ltimes G_1$.
\end{prop}
\begin{proof}
  The vector space $\quot(n)_{\sym_n}$ is one-dimensional for all $n$,
  with basis given by the image of the corolla with $n$ nodes. Let us
  denote this basis element by $\nu_{n-1}$. Any element of
  $\GG_{\quot}$ can be uniquely written as the product $\lambda
  (\sum_{m\geq 0} f_m \nu_m)$ of $\lambda \in \QQ^*$ and $f=\sum_{m\geq
    0} f_m \nu_m$ with $f_0=1$. Let us compute the product of $\lambda
  f=\lambda (\sum_{m\geq 0} f_m \nu_m)$ and $\theta g=\theta
  (\sum_{n\geq 0} g_n \nu_n)$ with the conventions $f_0=1$ and
  $g_0=1$. One finds that
  \begin{equation}
    \lambda f \times\theta g =
\sum_{m\geq 0}\sum_{n\geq 0} \lambda f_{m} \theta^{m}
    (\theta g_{n}) \nu_{n+m}= \lambda \theta \sum_{m\geq 0}  \sum_{n\geq
      0} \theta^{m}
    f_{m} g_{n} \nu_{n+m}.
  \end{equation}
  One defines a map from $\GG_{\quot}$ to $G_2 \ltimes G_1$ by
  $\lambda(\sum_m f_m \nu_m)\mapsto (\lambda,f(x))$ with $f(x)=\sum_m
  f_m x^m$. The product in $G_2 \ltimes G_1$ is given by
  \begin{equation}
    (\lambda, f(x))(\theta, g(x))=(\lambda \theta , f(\theta x) g(x)).
  \end{equation}
  Hence the map is an isomorphism.
\end{proof}

Let us denote by $\psi$ the quotient map from $\prelie$ to $\quot$.

\begin{prop}
  The image of $\exp^*$ by $\psi$ is $(\exp
  x-1)/x$ and the image of $\log^*$ is $x/(
\exp x-1)$.
\end{prop}
\begin{proof}
  One must prove that the coefficient of the corolla with $n$ nodes in
  $\exp^*$ is $1/n!$ for all $n$. The argument is a simple recursion
  using (\ref{recursive}).
\end{proof}

Therefore the coefficients of the corollas in $\log^*$ are related to
the Bernoulli numbers, whose exponential generating function is
precisely $x/( \exp x-1)$.

\section{Expansion}

The coefficients of $\exp^*$ are easily computed by the recursive
formula (\ref{recursive}). They are known as the Connes-Moscovici
coefficients, and there exists a direct procedure to compute the
coefficient of any rooted tree in $\exp^*$, see \cite[\S
2.2]{brouder}. It is an interesting open problem to compute the
coefficients of $\log^*$, for which there is no known alternative to
the inversion in the group $\GG_{\prelie}$. We give here the first
terms of the expansions of $\exp^*$ and $\log^*$ in the rooted tree
basis of $\widehat{\prelie}$.

\begin{multline}
\exp^*=\epsfig{file=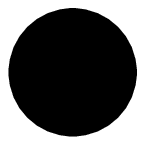,height=2mm}+
\frac{1}{2}\epsfig{file=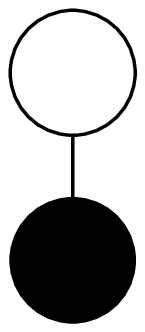,height=5mm}+
\frac{1}{6}\left(\epsfig{file=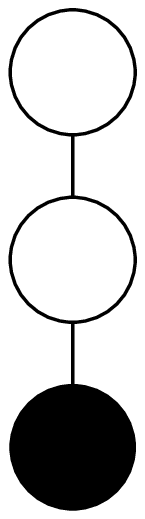,height=8mm}+\epsfig{file=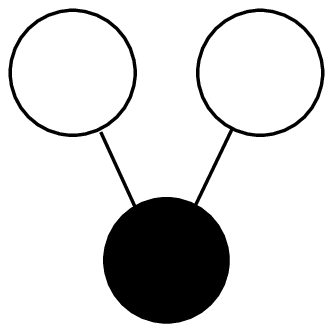,height=5mm}\right)+
\frac{1}{24}\left(\epsfig{file=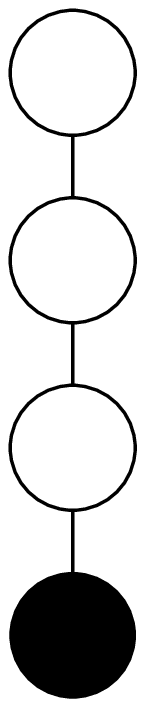,height=11mm}+\epsfig{file=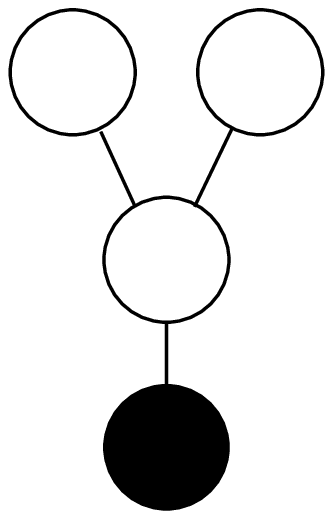,height=8mm}+3\epsfig{file=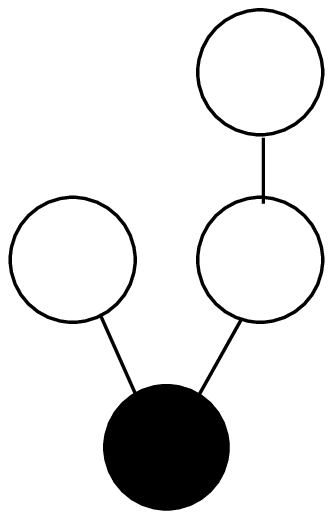,height=8mm}+\epsfig{file=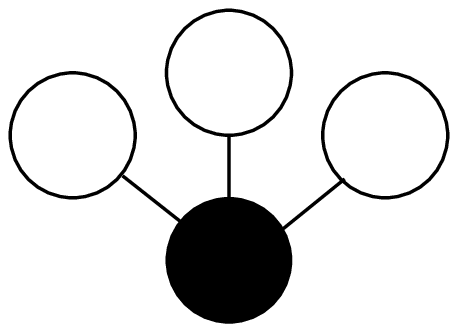,height=5mm}\right)+\\
\frac{1}{120}\left(\epsfig{file=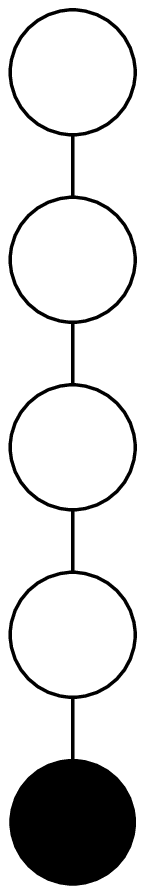,height=14mm}+\epsfig{file=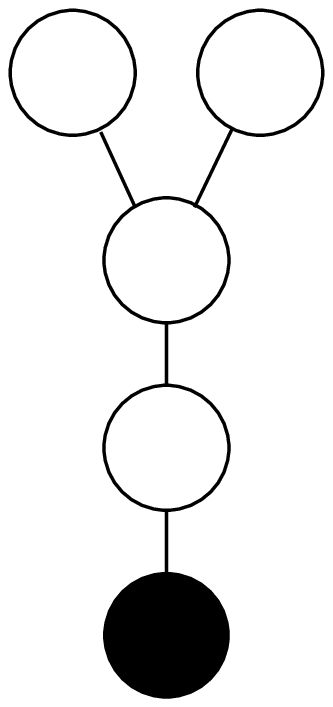,height=11mm}+3\epsfig{file=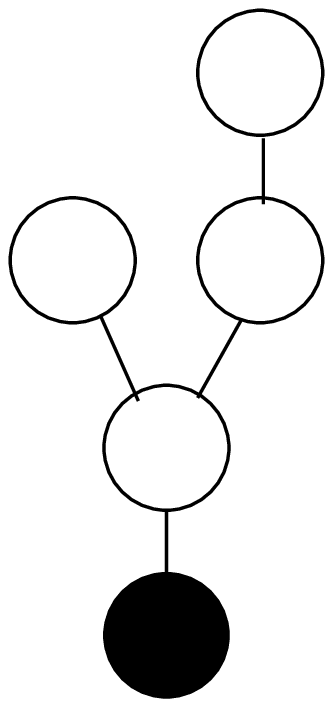,height=11mm}+\epsfig{file=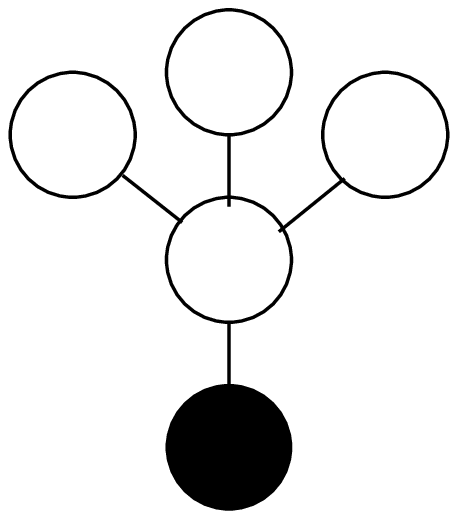,height=8mm}+3\epsfig{file=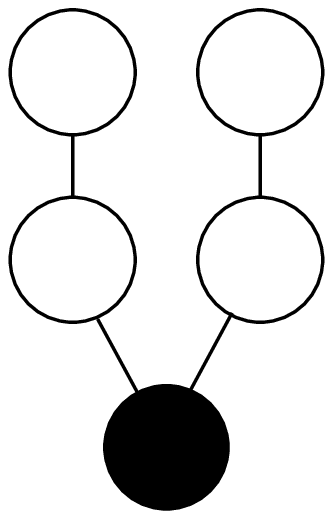,height=8mm}+4\epsfig{file=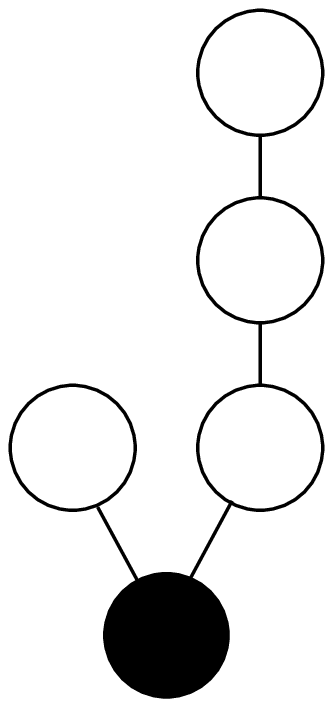,height=11mm}+4\epsfig{file=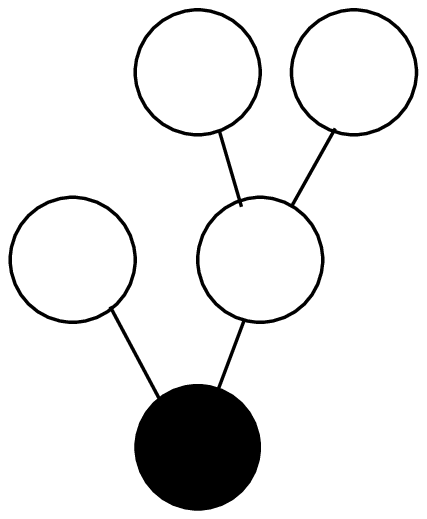,height=8mm}+6\epsfig{file=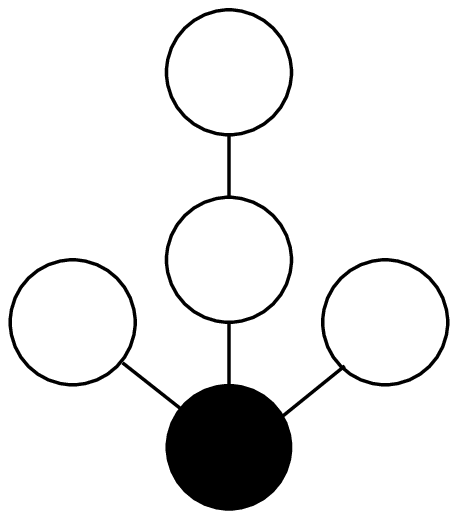,height=8mm}+\epsfig{file=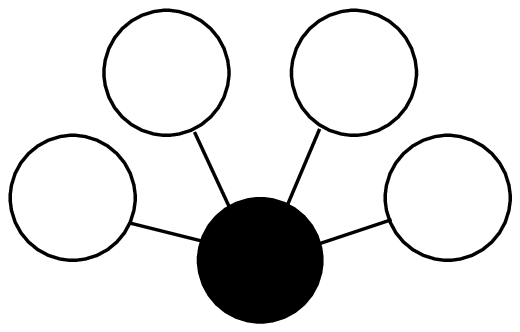,height=5mm}\right)+\,\cdots
\end{multline}

\begin{multline}
\log^*=\epsfig{file=a0.eps,height=2mm}-\frac{1}{2}\epsfig{file=a10.eps,height=5mm}+\frac{1}{6}\left(2
\epsfig{file=a110.eps,height=8mm}+\frac{1}{2}\epsfig{file=a200.eps,height=5mm}\right)-\frac{1}{24}\left(6
\epsfig{file=a1110.eps,height=11mm}+2\epsfig{file=a1200.eps,height=8mm}+2\epsfig{file=a2100.eps,height=8mm}\right)
+\\
\frac{1}{120}\left(24\epsfig{file=a11110.eps,height=14mm}+4\epsfig{file=a11200.eps,height=11mm}+12\epsfig{file=a12100.eps,height=11mm}+\frac{2}{3}\epsfig{file=a13000.eps,height=8mm}+2\epsfig{file=a21010.eps,height=8mm}+6\epsfig{file=a21100.eps,height=11mm}+\epsfig{file=a22000.eps,height=8mm}-1\epsfig{file=a31000.eps,height=8mm}-\frac{1}{30}\epsfig{file=a40000.eps,height=5mm}\right)
+\,\cdots
\end{multline}

\nocite{*}
\bibliographystyle{plain}
\bibliography{expon}

\end{document}